\documentclass{amsart}

\usepackage{amsfonts}
\usepackage{bezier}
\usepackage{euscript}
\usepackage{amsfonts,amssymb,amsmath,amsbsy,amsthm,amscd}
\usepackage{xcolor}
\usepackage{graphicx}
\usepackage{color}
\usepackage[T2A]{fontenc}
\usepackage[utf8]{inputenc}
\usepackage[english]{babel}   

\newtheorem{Theorem}{Theorem}
\newtheorem{cor}{Corollary}
\newtheorem{rem}{Remark}

\def\br{{\bf r}}

\def\dfrac#1#2{\displaystyle{#1\over #2}}

\def\bv{{\bf V}}
\def\bV{{\bf V}}

\def\Div{\mbox{div}\,}
\def\Rot{\mbox{curl}\,}

\def\bB{{\bf B}}
\def\bx{{\bf x}}
\def\bE{{\bf E}}
\def\bV{{\bf V}}
\begin{document}

\title[ 
Affine solutions of cold plasma equations]{
On the properties of affine solutions of cold plasma equations
}

\author{Olga S. Rozanova*}
\address{ Mathematics and Mechanics Department, Lomonosov Moscow State University, Leninskie Gory,
Moscow, 119991,
Russia}
\email{rozanova@mech.math.msu.su}

\author{ Marko K. Turzynsky}
\address{Russian University of Transport, Obraztsova, 9,   Moscow, 127055
and  Higher School of Economics, Pokrovskiy Blvd, 11, Moscow, 109028, Russia
}
\email{m13041@yandex.ru}





\maketitle


\begin{abstract}
We study the affine solutions of the equations of plane oscillations of cold plasma, which, under the assumption of electrostaticity, correspond to the Euler-Poisson equations in the repulsive case. It is proved that the zero equilibrium state of the cold plasma equations, both with and without the assumption of electrostaticity, is unstable  in the class of all affine solutions. It is also shown that an arbitrary perturbation of an axially symmetric electrostatic solution leads to a finite time blow-up.
\end{abstract}


\section{Introduction}

The equations that describe the cold plasma oscillations in the Euler coordinates have the form \cite{Textbook1}, \cite{Textbook2}:
\begin{equation}\label{1}
\partial_t  n +{\rm div}(n {\bf V})=0,\quad
\partial_t {\bf V} +  ({\bf V}\cdot \nabla ){\bf V}=-({\bf E}+[{\bf V} \times {\bf B}]),
 \end{equation}
 \begin{equation}\label{2}
\partial_t{\bf E} = n\!{\bf V}+{\Rot}\,{\bf B},\quad \partial_t{\bf B}=-{\Rot}\,{\bf E},\quad {\rm div}\,{\bf B}=0,
\end{equation}
where ${\bf V}(t,x)$, $n(t,x)>0$ are the speed and density of electrons,
${\bf B}(t, x)$ and ${\bf E}(t,x)$ are electrical and
magnetic fields, $x\in{\mathbb R}^3,$ $t\ge 0$, $\nabla$, $\rm div$, $\Rot$ are the gradient, divergence and vorticity with respect to the spatial variables.

System \eqref{1}, \eqref{2} has a class of solutions
\begin{equation}\label{VE}
\bV=Q(t) \bx, \qquad \bE=R(t)\bx,
 \end{equation}
 where $Q$ and $R$ are $3 \times 3$ matrices with coefficients depending on $t$, and $\bf x$ is the radius vector of the point $x\in \mathbb R^3$. These solutions are called {\it affine}.
Affine solutions have been known since the time of Kirchhoff and Dirichlet. They play an important role in various models of continuous media \cite{Borisov}.

It follows from  \eqref{VE} that $n=n(t)$ and $\bB=\bB(t)$. Since ${\Rot}\,{\bf B}(t)=0$, then from the first equations of  system \eqref{1}, \eqref{2}, under the assumption that the solution is sufficiently smooth and the steady-state density equals $1$, we get  $n=1-{\rm div}\,{\bf E}>0$. Thus, we can exclude $n$ from the system and  obtain
\begin{equation}\label{4}
\partial_t {\bf V} +  ({\bf V}\cdot \nabla ){\bf V}=-({\bf E}+[{\bf V} \times {\bf B}(t)]),\quad
\partial_t {\bf E} + {\bf V}{\rm div}\,{\bf E}={\bf V}, \quad \dot \bB (t)+{\Rot}\,{\bf E}=0.
\end{equation}

We will consider system \eqref{4}  together with the initial data
\begin{equation*}\label{CD}
({\bf V},\,{\bf E},\,{\bf B})|_{t=0}=(Q_0 \bx, R_0 \bx,{\bf B}_0),
\end{equation*}
with constant matrices $Q_0$ and $R_0$.

 Further, we assume that oscillations occur in a plane perpendicular to the coordinate vector $e_3$.
  Thus, we restrict the class of solutions under consideration to
 \begin{equation}\label{6}
  {\bf V}=Q \bx=\left(\begin{array}{ccr} a(t) & b(t) & 0\\
c(t) & d(t) & 0 \\ 0 & 0 & 0
\end{array}\right) \bx,\quad {\bf E}=R \bx=\left(\begin{array}{ccr} A(t) & B(t) & 0\\
C(t) & D(t) & 0 \\ 0 & 0 & 0
\end{array}\right) \bx,\quad \bB=(0, 0, \mathcal B(t)).
 \end{equation}

System \eqref{4} for solutions of the form \eqref{6} reduces to a matrix system of differential equations for the matrices $Q$ and $R$ and the scalar function $\mathcal B$:
\begin{equation}\label{7}
 \dot Q + Q^2-\mathcal B(t) \mathcal{L}Q+R=0,\quad \dot R-(1-{\rm tr}\,R)Q=0,\quad  \dot {\mathcal B}(t)-{\rm tr}\,(\mathcal{L}R)=0,
\end{equation}
which consists of 9 differential
equations, here $\mathcal{L}=\left(\begin{array}{ccr} 0 & -1 & 0\\
1 & 0 & 0 \\ 0 & 0 & 0
\end{array}\right)$.
The initial data for \eqref{7}  are
\begin{equation*}\label{CD}
(Q,\,R,\,{\mathcal B})|_{t=0}=(Q_0 , R_0 ,{\mathcal B}_0).
\end{equation*}

An important class of oscillations is distinguished by the condition ${\Rot}\,{\bf E}=0$, such oscillations are called {\it electrostatic}. Under this condition, the magnetic field does not change with time and has the form
  ${\bf B}=(0, 0, \mathcal B_0)$. Another consequence of this assumption is the condition ${\Rot}\,{\bf V}=0$.
Thus, system \eqref{4} can be rewritten as
\begin{equation}\label{3}
 \partial_t  n +{\rm div}(n {\bf V})=0,\quad
\partial_t {\bf V} +  ({\bf V}\cdot \nabla ){\bf V}=-({\bf E}+[{\bf V} \times {\bf B}_0]),\quad
\partial_t{\bf E} = n{\bf V},\quad {\Rot}\,{\bf E}=0.
\end{equation}
It is easy to check that this situation is realized only in the case $\mathcal B_0=0$, $b=c=B=C=0$, the number of equations in system \eqref{7} is reduced to four.

If we introduce a potential $\Phi$ such that $\nabla \Phi = -\bE$, then  \eqref{3} can be rewritten as a system of Euler-Poisson equations
(e.g. \cite{ELT})
\begin{eqnarray}\label{EP}
\dfrac{\partial n }{\partial t} + \Div(n \bv)=0,\quad
\dfrac{\partial \bv }{\partial t} + \left( \bv \cdot \nabla \right)
\bv =\,k \,  \nabla \Phi, \quad \Delta \Phi =n-n_0,
\end{eqnarray}
 for $k=n_0=1$. Thus, the results obtained for system \eqref{3}  in the electrostatic case can be reformulated in terms of solutions of the Euler-Poisson equations.

 Solutions of the form \eqref{6}  have also a subclass of solutions  {\it with the radial symmetry in the plane $x_3=0$}, for which
\begin{equation}\label{symm}
 \bV=a \,\br + c\, \br_\bot, \quad \bE= A\, \br+ C\,\br_\bot, \quad \br= (x_1, x_2, 0), \quad \br_\bot= (x_2, -x_1, 0).
\end{equation}
For such solutions, the number of equations in  system \eqref{7}  is reduced to five. Such solutions are electrostatic only if $c=C={\mathcal B}=0$, i.e. the radially symmetric electrostatic solution  is axisymmetric.

It was recently proved \cite{R22_Rad} that affine solutions play an exceptional role in the class of axisymmetric solutions of multidimensional Euler-Poisson equations \eqref{EP}  depending on $(t, r)$, where $r=\sqrt{x_1^2+x_2^2}$. Namely, if some solution preserves global smoothness in time, then it is either affine or tends to affine as $t\to \infty$ uniformly on each interval in $r$. In addition, the zero equilibrium state turns out to be unstable with respect to axisymmetric perturbations of an arbitrary form, but stable with respect to affine axisymmetric perturbations. As shown above, the axisymmetric solutions of the Euler-Poisson equations correspond to the axisymmetric solutions of the cold plasma equations with the condition of electrostaticity. In this regard, a natural question arises: will the zero equilibrium of the cold plasma equations be stable in the class of affine solutions without the assumption of axial symmetry or electrostaticity?

We show that the answer to this question is negative. Moreover, it turned out that a general perturbation from the affine axially symmetric solution leads to a blow-up of the solution in a finite time. Since plane oscillations are a subclass of spatial oscillations, the result on the instability of the zero equilibrium state is also valid for the three-dimensional case.

The paper has the following structure. In Section \ref{S2} for the electrostatic case $\mathcal B_0=0$ we construct a globally smooth solution of system \eqref{8} with axial symmetry and show that it is Lyapunov stable in the class of electrostatic solutions \eqref{6} with axial symmetry. In particular, the equilibrium $\bV=\bE=0$ of system \eqref{8}  is stable. Further, relying on the Floquet theory, we  show that the equilibrium $\bV=\bE=0$ is unstable with respect to small perturbations in the class of affine electrostatic solutions, and thus also in the class of arbitrary affine perturbations. Moreover, we  show that in the general case (without a special choice of initial data) such a perturbation grows with time and leads to a blow-up of the solution in a finite time. Further, by numerical computation of the characteristic multipliers of the system of ordinary differential equations, we show that a similar result is valid for an arbitrary deviation from the globally smooth solution with axial symmetry constructed in Section \ref{S2}.

In Section \ref{S3} we consider the non-electrostatic case and show that the zero equilibrium $\bV=\bE=\mathcal B_0=0$ is unstable in the class of radially symmetric non-electrostatic solutions. We  show that a general perturbation of a globally smooth axisymmetric electrostatic solution  in the class of radially symmetric non-electrostatic solutions also leads to a blow-up of the solution in a finite time.
Also in section \ref{S3} we  discuss the difference between deviations from zero equilibrium $\bV=\bE=\mathcal B_0=0$ and equilibrium $\bV=\bE=0$, $\mathcal B_0\ne 0 $.

\section{Electrostatic oscillations} \label{S2}

\subsection  {Solution with axial symmetry}\label{S2.1}
In the case  \eqref{symm}, under the electrostatic condition $c=C=\mathcal B_0=0$, the system \eqref{7} takes the form
\begin{equation}\label{8}
  \dot a=-A-a^2,\quad  \dot A=a-2Aa.
\end{equation}
The equilibrium of the $a=A=0$ system is the center, since the eigenvalues of the linear approximation matrix in it are equal to $\lambda_{1,2}=\pm i$, and the phase curves are symmetric when $a$ is replaced by $-a $.

 Let us find the first integrals of  system \eqref{8}. It implies
$ \frac{a da}{dA}-\frac{a^2}{2A-1}=\frac{A}{2A-1},$
after replacing $a^2=u$ we obtain
a linear equation, the solution is
 \begin{equation}\label{91}
 a=\pm\sqrt{(\tfrac{1}{2}\ln|2A-1|+K)(2A-1)-\tfrac{1}{2}}, \quad K=\rm const.
 \end{equation}

This integral was also obtained in \cite{R22_Rad}. The curve on the phase plane given by the expression \eqref{91} is bounded for all values of $K$ satisfying the condition $A(0)<\frac12$ (see \cite{R22_Rad}), so the derivatives of solution, $a$ and $A$, remain bounded for all $t>0$. The explicit form of the bounded phase curve  \eqref{91} also implies that the equilibrium $a=A=0$, corresponding to the zero rest state, is stable by Lyapunov. The solutions corresponding to any fixed phase curve \eqref{91} are also Lyapunov stable if the perturbation occurs in the class of affine solutions given by  system \eqref{8}.

In this way, $a(t)$, $A(t)$ are periodic with period
\begin{equation*}\label{T}
  T= 2 \int\limits_{A_-}^{A_+} \frac {d\eta}{(1-2 \eta )a(\eta)},
\end{equation*}
$a$ is set to \eqref{91}, $A_-<0$ and
$A_+>0$ is the smaller and larger roots of the equation $a(A)=0$.
Besides, $\int\limits_0^{T} a(\tau) \, d\tau=0$.

The period $T$ was studied in \cite{R22_Rad}. It depends on $A(0)=\varepsilon$, $\varepsilon\in (0, \frac12)$ decreasing monotonically from $2\pi$ to $\sqrt{2}\pi$, and the asymptotic formula
 \begin{equation}\label{T_per}
T=2\pi (1- \frac{1}{12} \varepsilon^2+o(\varepsilon^2)),\quad \varepsilon\to 0,
\end{equation}
holds, see
 \cite{R22_Rad}, Lemma 4.

\subsection  {Arbitrary electrostatic oscillations of  form \eqref{6}}

We formulate two similar theorems, the first of which will be proved analytically, while the second is a semi-analytical result, for the proof we use numerical methods.

The proof of all theorems is based on the Floquet theory for systems of linear equations with periodic coefficients (for example, \cite{Chicone}, Section 2.4). According to this theory, for the fundamental matrix $\Psi (t)$ ($\Psi (0)=E$) there exists a constant matrix $M$, possibly with complex coefficients, such that $\Psi (T)=e^{T M}$, where $T$ is the period of the coefficients. The eigenvalues of the matrix of monodromy $e^{T M}$ are called the characteristic multipliers of the system. If among the characteristic multipliers there are such that their absolute value is greater than one, then the zero solution  of the studied linear system is unstable in the sense of Lyapunov (\cite{Chicone}, Theorem 2.53).

\begin{Theorem}\label{T1}
1. The zero equilibrium of  system \eqref{7} in the class ${\mathcal B}(t)\equiv 0$ (corresponding to electrostatic oscillations) is unstable by Lyapunov.

2. A general  small non-axisymmetric perturbation of the equilibrium blows up in a finite time.
\end{Theorem}

{\it Proof of Theorem \ref{T1}.}
The system \eqref{7} in the case of ${\mathcal B}\equiv 0$ has the form
\begin{equation}\label{sym0}
 \dot A=(1-A-D)a,\quad \dot D=(1-A-D)d,\,\quad \dot a+a^2+A=0,\quad \dot d+d^2+D=0.
\end{equation}
To study the effect of deviation from symmetry, we make the substitution $d=a+\sigma$, $D=A+\delta$, which corresponds to the axisymmetric case for $\sigma=\delta=0$. Then \eqref{sym0} reduces to
\begin{equation*}\label{9}
 \dot A=(1-2A)a-\delta a,\quad \dot a=-a^2-A,\,\quad \dot \delta=(1-2A-\delta)\sigma,\quad \dot \sigma=-\sigma^2-2a\sigma-\delta.
\end{equation*}
We choose a small parameter $\varepsilon$ and set
\begin{eqnarray*}\label{101}
  A(t)&=&  A_0(t)+\varepsilon^2  A_1(t)+o(\varepsilon^2),\quad a(t)=  a_0(t)+\varepsilon^2 a_1(t)+o(\varepsilon^2), \\ \delta(t)&=&\varepsilon^2 \delta_1(t)+o(\varepsilon^2), \quad\sigma(t)=\varepsilon^2 \sigma_1(t)+o(\varepsilon^2),\nonumber
\end{eqnarray*}
For $\varepsilon=0$ we obtain a globally smooth solution $A_0(t), a_0(t)$, which is a solution to  system \eqref{8}.
 For the functions $A_1,\, a_1,\,\delta_1,\, \sigma_1$, discarding terms of the order of smallness $o(\varepsilon^2)$, we obtain the linear system
\begin{eqnarray}\label{s4.1}
\dot A_1&=&- 2 a_0 A_1 + (1-2A_0)a_1-a_0 \delta_1, \quad \dot a_1=-2a_0 a_1-A_1,\\
\dot \delta_1&=&(1-2A_0)\sigma_1, \quad \dot \sigma_1=-2a_0\sigma_1-\delta_1.\label{s4.2}
\end{eqnarray}

Let us show that the zero solution of  system \eqref{s4.1}, \eqref{s4.2} is unstable.

 Note that if $\delta_1(0)=\sigma_1(0)=0$, then  system \eqref{s4.1}, \eqref{s4.2} reduces to two equations \eqref{s4.1}, $ \delta_1\equiv 0$, and then the equilibrium $A_1=a_1=0$ turns out to be stable. This follows from the fact that the perturbed solution remains axisymmetric and the integral \eqref{91} holds for it.
 Thus, we will consider a perturbation of the solution $A_0(t), a_0(t)$, for which $\delta_1(0)\sigma_1(0)\ne 0$.

 1. Let $A_0(t), a_0(t)$ itself be a small deviation from the zero equilibrium position. We take $A_0(0)=\varepsilon$ as a small parameter. Then
  \begin{eqnarray}\label{Aa0}
  \qquad A_0(t)=\varepsilon \cos t + A_{01}(t)\varepsilon^2 +o(\varepsilon^2), \quad a_0(t)=-\varepsilon \sin t + a_{01}(t)\varepsilon^2 +o(\varepsilon^2),
\end{eqnarray}
all subsequent expansion terms are found sequentially from  \eqref{8}.

To obtain an asymptotic representation of the components of the fundamental matrix, we set
 \begin{eqnarray}\label{A1}
  A_1(t)&= &A_{10}(t)+A_{11}(t)\varepsilon+A_{12}(t)\varepsilon^2+o(\varepsilon^2),\\
   a_1(t)&=& a_{10}(t)+a_{11}(t)\varepsilon+a_{12}(t)\varepsilon^2+o(\varepsilon^2),\label{a1} \nonumber\\
    \delta_1(t)&=&  \delta_{10}(t)+ \delta_{11}(t)\varepsilon+ \delta_{12}(t)\varepsilon^2+o(\varepsilon^2),\label{d1} \nonumber\\
    \sigma_1(t)&=& \sigma_{10}(t)+\sigma_{11}(t)\varepsilon+\sigma_{12}(t)\varepsilon^2+o(\varepsilon^2).\label{s1}
\end{eqnarray}
The fundamental matrix has the form $\Psi(t)=\Psi_0(t)+ \Psi_1(t) \varepsilon + \Psi_2(t) \varepsilon^2+ o(\varepsilon^2),$
$\Psi_0(0)=\mathbb E$, $\Psi_i(0)=0$, $i\in \mathbb N$. The calculations that need to be done to find the matrices $\Psi_i(t)$  are cumbersome, but  standard: we substitute \eqref{Aa0}, \eqref{A1} -- \eqref{s1} into \eqref{s4.1}, \eqref{s4.2} and equate the coefficients at the same powers of $\varepsilon$. At each stage, one has to solve a linear inhomogeneous system with constant coefficients.

Denote the eigenvalues of the matrix $\Psi(T)$ as $\lambda_i,$ and the eigenvalues of the matrix $\Psi_k(T)=\Psi_0(T_0)+ \Psi_1(T_1) \varepsilon +\cdots +\Psi_k(T_k) \varepsilon^k$, $k\in \mathbb N$ as $\bar \lambda_{ki},$ $\, i=1,\dots, 4$. Further, we denote  as $T_j$, $j=0, \dots, k$, the period $T$ (see \eqref{T_per}) calculated in the approximation $O(\varepsilon^j)$, and the eigenvalues of the matrix $\Psi_j(T_j)$ as $\lambda_{ji }$.
   To prove the instability, we have to find an expansion up to such an order  $k$ in $\varepsilon$ that among $\bar \lambda_{ki}$ there is a greater than one by the absolute value.

 As follows from \eqref{T_per}, $T_0=T_1=2 \pi$, $T_2=2 \pi-\frac{\pi}{6} \varepsilon^2$.
 Calculations performed using the computer algebra package MAPLE  show  that $\bar \lambda_{0i}=\bar \lambda_{1i}=1,$  $i=1,\dots, 4$,
 \begin{eqnarray}\label{EV4} \qquad
 \bar \lambda_{2i}= 1\pm\frac{\sqrt{3} \pi}{6}\varepsilon^2+ o (\varepsilon^2),\, i=1,2, \quad \bar \lambda_{2i}= 1 \pm \frac{\sqrt{3} \pi}{2}\varepsilon^2 + o (\varepsilon^2),\, i=3,4,
\end{eqnarray}
and  further  terms of the expansion cannot change the coefficients of $\varepsilon$ to a power less than two. Thus, already for $k=2$ we can conclude that there is a pair of eigenvalues such that $| \lambda|>1$. Thus,  the instability of the zero equilibrium is proved.

Note that according to the Liouville theorem on the conservation of the phase volume
 \begin{eqnarray*}
 {\rm det} \Psi(T)= \exp \left(\int\limits_0^T {\rm tr} {\mathcal M} (\tau) d\tau \right)\,{\rm det} \Psi(0),
 \end{eqnarray*}
where $\mathcal M$ is the matrix corresponding to the linear system
\eqref{s4.1}, \eqref{s4.2}. It is easy to see that ${\rm tr} {\mathcal M}= 6 a_0(t)$, and since  $\int\limits_0^T a_0 (\tau) d\tau=0$. Therefore $\prod\limits_{i=1}^4 \, |\lambda_i|\,=1.$ However
$\prod\limits_{i=1}^4 \, |\bar\lambda_{ki}|$ must be equal to one only up to terms  $o(\varepsilon^k)$, which we see from \eqref{EV4}.

2.  The fact that the linear system \eqref{s4.1}, \eqref{s4.2} has at least one characteristic multiplier whose absolute value is greater than one indicates that any component of its solution, including $A_1 (t)$, contains a term $  \mathcal C P(t) \exp(\mu t)$, with characteristic exponent $\mu>0$, $P(t)$ is a bounded periodic function, $ \mathcal C$ is a constant , depending on the initial data.
With some special choice of initial data, one can make the constant $ \mathcal C$ equal to zero. This will certainly be the case if $\sigma(0)=\delta(0)=0$, that is, when the initial perturbation is axisymmetric. However, if the perturbation is chosen arbitrarily, then an exponentially growing term is necessarily present. Therefore, if we assume the boundedness of $A(t)$, and hence $A_1(t)$, for all finite times $t$, then there exists a point $t_*$ for which $A(t_*)= \frac12$. Since $A(t)=\frac12$ is a part of solution to  system \eqref{7}, this would contradict the  uniqueness theorem. Therefore, there is a finite time $t_c<t_*$ at which $A(t)$ becomes infinite, which corresponds to a blow-up of the solution.  $\square$

\begin{Theorem}\label{T2}
A globally smooth axisymmetric solution of system \eqref{7}  is  unstable by Lyapunov in the class ${\mathcal B}(t)\equiv 0$, and any general non-axisymmetric perturbation of it blows up in a finite time.
\end{Theorem}

{\it Proof of Theorem \ref{T2}} is completely analogous to the proof of the previous theorem, but in order to investigate the stability of the perturbation from the equilibrium  of the solution $A_0(t), a_0(t)$, one has to apply numerical methods.
This is not surprising, since even for a much simpler situation the Mathieu equations  the stability region can only be found numerically without the assumption that the periodic coefficient is small \cite{Mathieu}.
Therefore, for each fixed $A_0(0)=A_*$, we solve  system \eqref{8}, \eqref{s4.1}, \eqref{s4.2} numerically using the Runge–Kutta–Felberg method of the fourth-fifth order (RKF45),  and then find the absolute values of the eigenvalues at the point $T(A_*)$.

\begin{figure}[htb!]
\hspace{-1.0cm}
\begin{minipage}{0.3\columnwidth}
\includegraphics[scale=0.3]{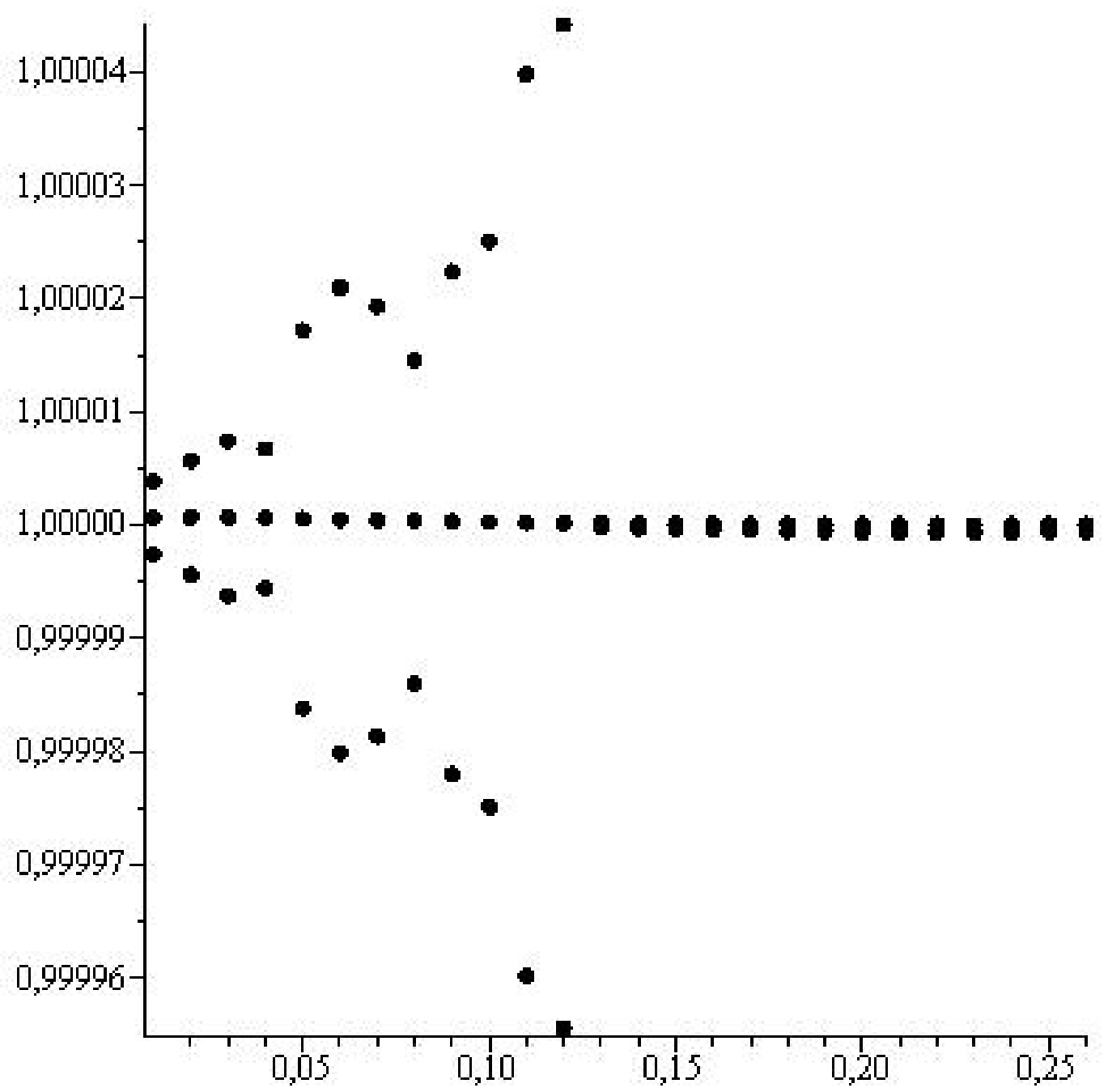}
\end{minipage}
\hspace{0.5cm}
\begin{minipage}{0.3\columnwidth} 
\includegraphics[scale=0.6]{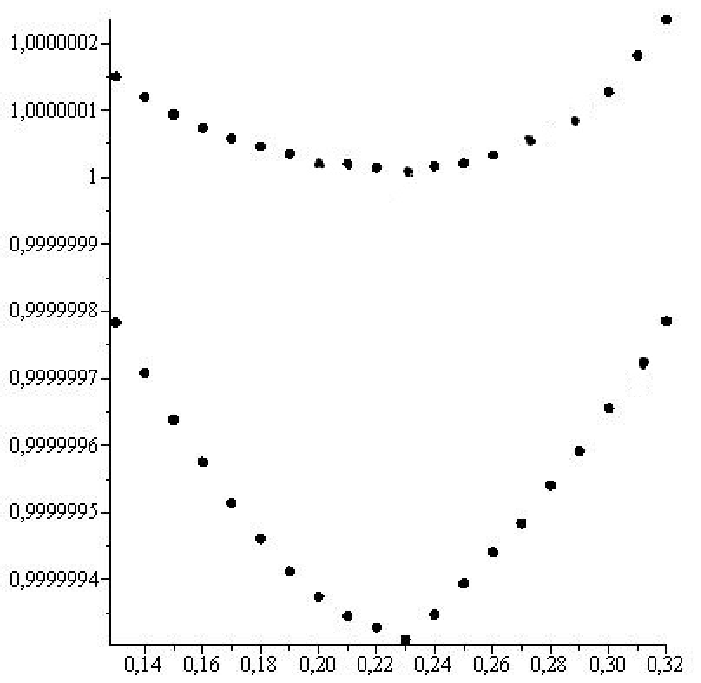}
\end{minipage}
\hspace{0.5cm}
\begin{minipage}{0.3\columnwidth} 
\includegraphics[scale=0.3]{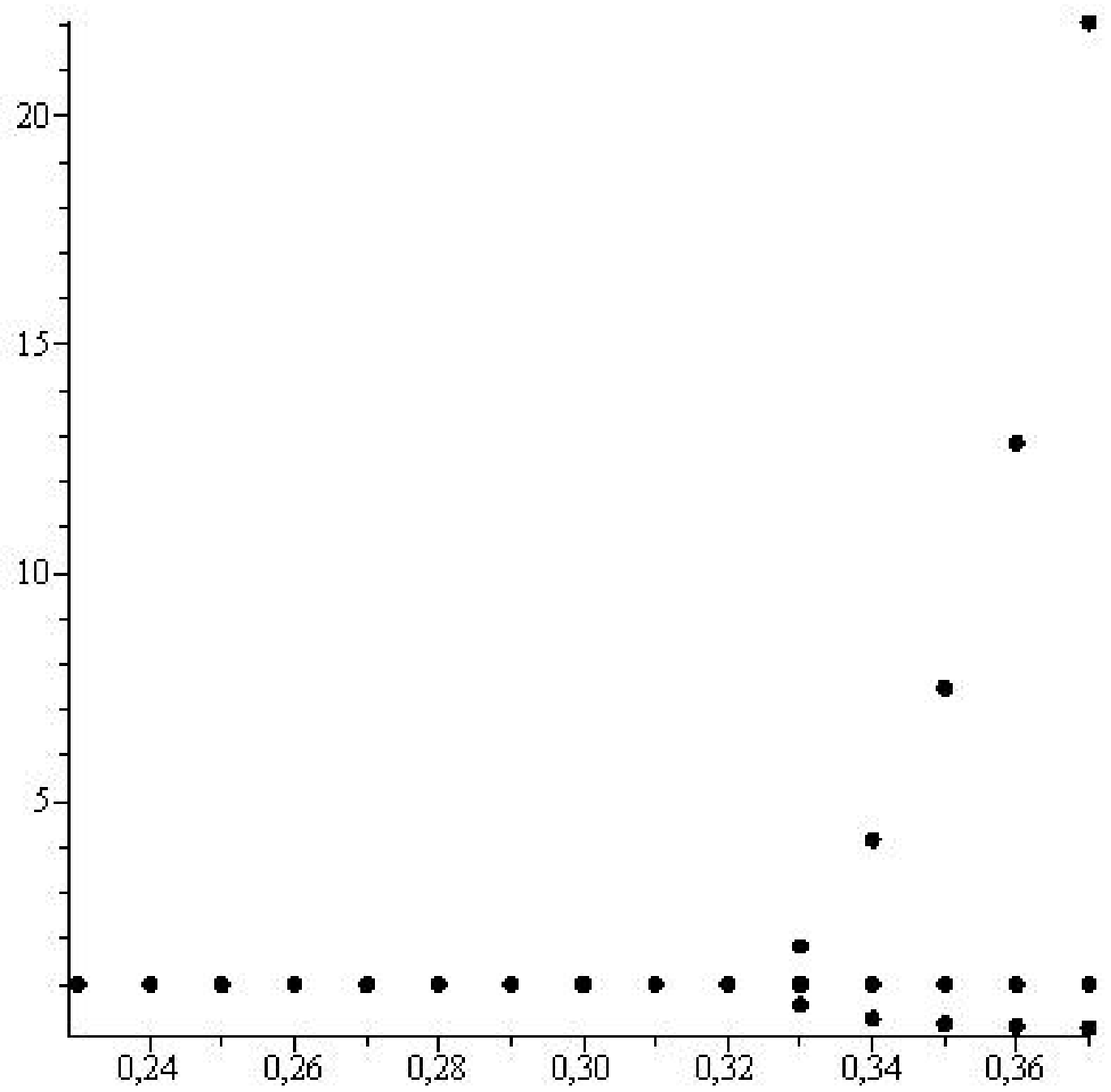}
\end{minipage}
\caption{Values of characteristic multipliers for the case of electrostatic non-axisymmetric oscillations depending on $A(0)$. The figure in the center shows that the curves that appear to match in the figure on the left are actually different.}
\label{Pic1}
\end{figure}

 Figure \ref{Pic1} shows the dependence  $|\lambda_i(A_*)|$ for different ranges of $A_*$. It is easy to see that the maximum of the eigenvalues is greater than 1 in absolute values, which indicates instability.

If we introduce the measure of instability as $S(A_*)=\max\limits_i |\lambda_i(A_*)| -1$, then we see that this value on the interval $(0, \frac12)$ varies nonmonotonically, remaining positive. Namely, it initially increases for small $A_*$, which is indicated by the asymptotic representation of the eigenvalues \eqref{EV4}, but sharply decreases at the point $A_*\approx 0.125$, remaining very small up to the point $A_*\approx 0.32$, after which it sharply increases, approaching the boundary value $A_*\approx 0.5$. This, in particular, indicates that it is very difficult to detect instability by direct numerical methods without resorting to the Floquet theory. since the deviation from the stationary solution $A_0(t), a_0(t)$ grows very slowly. Indeed, for example, in the region of the point $A_*=0.25$ we have $S(A_*)\sim 10^{-7}$.

We also note the properties of the eigenvalues themselves. On the interval $A_*\in(0,\approx 0.125)$ there are a pair of complex conjugate and a pair of real eigenvalues, one of which is greater than one. On the interval $A_*\in(\approx 0.125, \approx 0.32)$, there are two pairs of complex conjugate eigenvalues, which also indicates a softer loss of stability. On the interval $A_*\in(\approx 0.32, 0.5)$, a pair of real eigenvalues again arises, one of which is positive and rapidly grows as one approaches the right boundary of the interval.

The blow-up of an arbitrary non-axisymmetric perturbation of a globally smooth axisymmetric solution is proved in the same way as in Theorem \ref{T1}.
$\square$

\bigskip
\begin{cor}\label{c1}  The zero equilibrium of  system \eqref{4}  is unstable by Lyapunov  in the class of all  affine solutions \eqref{VE}.
\end{cor}

{\it The proof} of Corollary \ref{c1} follows directly from the observation that plane deviations from the equilibrium  in the class  \eqref{6} with  ${\mathcal B}(t)\equiv 0$ form a subclass among all possible affine deviations. $\square$

 \begin{cor}\label{c2}  1. The zero equilibrium state of the Euler-Poisson equations is unstable in the class of affine solutions.

  2. A general small non-axisymmetric affine perturbation of a globally smooth axisymmetric affine solution of the Euler-Poisson equations \eqref{EP} blows up in a finite time.
\end{cor}

{\it Proof.}
Сorollary \ref{c2} is a reformulation of Theorem \ref{T1} for the case of the Euler-Poisson equations.

\begin{rem} The electrostatic solutions of the cold plasma equations, generally speaking, blow up in a finite time,  this time can be estimated from below, see \cite{R2}.
\end{rem}

\section{Non-electrostatic oscillations} \label{S3}

The equilibrium  of system \eqref{7}  corresponding to nonzero density have the form $Q=R=0$ (a matrix with zero components), $\bB=\bB_0=\rm const$.

The linearization matrix at this equilibrium  has the following  eigenvalues
$\lambda=\pm\frac{1}{2}
\sqrt{-4-2\mathcal B_0^2\pm2\sqrt{\mathcal B_0^4+4  \mathcal B_0^2}}, $
double multiplicity,
and $\lambda=0$.
 It is easy to check that $(-4-2 \mathcal B_0^2+2\sqrt{\mathcal B_0^4+4 \mathcal B_0^2})<0$ for all real $\mathcal B_0$.
Hence,  the real part of all eigenvalues is zero, and the theory of linear approximation to study the stability of equilibrium is not applicable.

Since we are interested in the deviation from the electrostatic condition when $\mathcal B_0=0$, we will study the stability of the equilibrium with $\mathcal B_0=0$ in the class of non-electrostatic perturbations.


\begin{Theorem}\label{T3}
1. The zero equilibrium of  system \eqref{7} is unstable in the sense of Lyapunov in the class of affine non-electrostatic solutions.

2. A general  small radially symmetric affine non-electrostatic perturbation of a globally smooth axisymmetric affine solution of  system \eqref{7} blows up in a finite time.
\end{Theorem}

\bigskip

{\it Proof of Theorem 3.} 1. The proof is completely analogous to the proof of Theorem 1 and is based on the Floquet theory described above. It suffices to show that the zero equilibrium is unstable in the class of affine solutions with radial symmetry \eqref{symm}.

System \eqref{7} in this case has the form
\begin{eqnarray}
\dot A &-& ( 1-2\,A ) a  =0,  \quad \dot C  -
 ( 1-2\,A ) c =0, \quad \dot {\mathcal B}  -2\,C  =0,\label{7.1}\\
\dot a&+&a^{2} -
c ^2   +A -\mathcal B c=0, \quad \dot c + 2\,c a  +C  +\mathcal B a  =0.\label{7.2}
\end{eqnarray}

Let us set
\begin{eqnarray}\label{10}
 A(t)&=&A_0(t)+ A_1(t)\varepsilon^2+o(\varepsilon),\quad a(t)=a_0(t)+ a_1(t)\varepsilon^2+o(\varepsilon^2), \\ c(t)&=& c_1(t)\varepsilon^2+o(\varepsilon^2), \quad C(t)=C_1(t)\varepsilon^2+o(\varepsilon^2),\nonumber \\ \mathcal B(t)&=&\mathcal B_1(t)\varepsilon^2+o(\varepsilon^2),\nonumber
\end{eqnarray}
where $\varepsilon$ is some small parameter. Substituting these series into \eqref{7.1}, \eqref{7.2} and remembering that \eqref{symm} implies $A=D$, $a=d$, $C=-B$, $c=-b$, we get the following system:
\begin{eqnarray}
 \dot a_0 &+&A_0 + a_0^2=0, \quad \quad \,\, \dot A_0 - ( 1-2 A_0) a_0=0,\nonumber\\
\dot a_1&+&
A_1 +2 a_0 a_1 =0, \quad
\dot A_1 - ( 1-2 A_0) a_1 +2  a_0 A_1 =0,\label{Aa11}\\
\dot C_1 & -& ( 1-2 A_0 )  c_1 =0,\quad
 \dot c_1  +a_0 {\mathcal B}_1 +2 a_0 c_1 +C_1 =0,\quad
\dot {\mathcal B}_1 -2 C_1  =0.\label{Cc11}
\end{eqnarray}
We see that the zero term of the series, $A_0(t), \,a_0(t)$, is a solution of \eqref{8}. For the next  terms of expansion, we obtain a linear system of equations \eqref{Aa11}, \eqref{Cc11} with known periodic coefficients. The equations \eqref{Aa11} for $A_1(t), \,a_1(t)$
is split off and  three equations \eqref{Cc11} can be considered separately.

  We choose $A_0(0)=\varepsilon\ll 1$, $a_0(0)=0$, so the zero terms of the series themselves turn out to be small,
  and the expansion \eqref {Aa0} is valid.

To obtain an asymptotic representation of the components of the fundamental matrix, we set
 \begin{eqnarray*}
   c_1(t)&=& c_{10}(t)+c_{11}(t)\varepsilon+c_{12}(t)\varepsilon^2+o(\varepsilon^2), \\  C_1(t)&=& C_{10}(t)+C_{11}(t)\varepsilon+C_{12}(t)\varepsilon^2+o(\varepsilon^2), \\ B_1(t)&=& B_{10}(t)+B_{11}(t)\varepsilon+B_{12}(t)\varepsilon^2+o(\varepsilon^2).\nonumber
\end{eqnarray*}
We use the same notation and methods as in the proof of Theorem \ref{T1}. Computations show that
  $\bar \lambda_{0i}=\bar \lambda_{1i}=1,$ $i=1, 2, 3$,
 \begin{eqnarray*}\label{EV3}
 \bar \lambda_{2i}=1\pm \frac{\sqrt{5}}{3}\pi\varepsilon^2+ o (\varepsilon^2),\, i=1,2,  \quad \bar \lambda_{23}=1,
  \end{eqnarray*}
and  further  terms of the expansion cannot change the coefficients of $\varepsilon$ to a power less than two.
Thus, there is a pair of eigenvalues such that $| \lambda|>1$ and the instability of the zero equilibrium is proved.

2. In order to prove the blow-up of an arbitrary non-electrostatic perturbation of a globally smooth axisymmetric solution, we cannot directly use the arguments of Theorem \ref{T1}. Indeed, our conclusions concern the components $C, c, \mathcal B$, while the restriction on the component $A$ led to a contradiction. Therefore, we note that the terms of expansion \eqref{10} for $a$ and $A$ in $\varepsilon$ starting from the fourth power, that is, $A_3$ and $a_3$, are no longer separated from $C, c, \mathcal B$. Namely, as follows from \eqref{7.1}, \eqref{7.2},  they are subject to the following inhomogeneous system of linear equations
\begin{equation*}
\dot A_3-a_3=-2 A_0 a_2-2 A_1 a_1-2 a_0 A_2,\quad
\dot a_3 +A_4= -2 a_0 a_2+a_1^2+{\mathcal B}_1 c_1-c_1^2.
\end{equation*}
Moreover, $A_0, A_1, A_2$, $a_0, a_1, a_2$, are periodic and bounded (which follows from \eqref{91}, since only the previous components $a$ and $A$ are used to calculate the first expansion components), while ${\mathcal B}_1$ and $ c_1$ generally contain the term $ \mathcal C P(t) \exp(\mu t)$, with characteristic exponent $\mu>0$, where $P(t)$ is a bounded periodic function, $ \mathcal C$ is a constant depending on the initial data.
With some special choice of initial data, it is possible to ensure that the constant $ \mathcal C$ turns out to be zero, but for an arbitrary non-electrostatic perturbation of the zero state of rest, an exponentially growing component of the solution is necessarily present.

 Therefore, as follows from the standard formula for representing the solution of a linear inhomogeneous equation, $A_3$ and $a_3$ also have this property, so if we assume that the solution is defined for all $t>0$, we get a contradiction with the condition $A< \frac12$.
$\square$

\bigskip
\begin{Theorem}\label{T4}
A globally smooth axisymmetric solution of  system \eqref{7} is unstable in the sense of Lyapunov in the class of affine non-electrostatic solutions with radial symmetry \eqref{symm}, and any general radially symmetric non-electrostatic perturbation of it  blows up in a finite time.
\end{Theorem}

For {\it the proof}, we repeat the procedure similar to  the proof of Theorem \ref{T2}.
For each fixed $A_0(0)=A_*$, we  solve  system \eqref{8}, \eqref{Cc11} numerically using the fourth-fifth order Runge–Kutta–Felberg method (RKF45), and then find the absolute values of eigenvalues at the point $T(A_*)$.

\begin{figure}[htb!]
\hspace{1.0cm}
\begin{minipage}{0.45\columnwidth}
\includegraphics[scale=0.3]{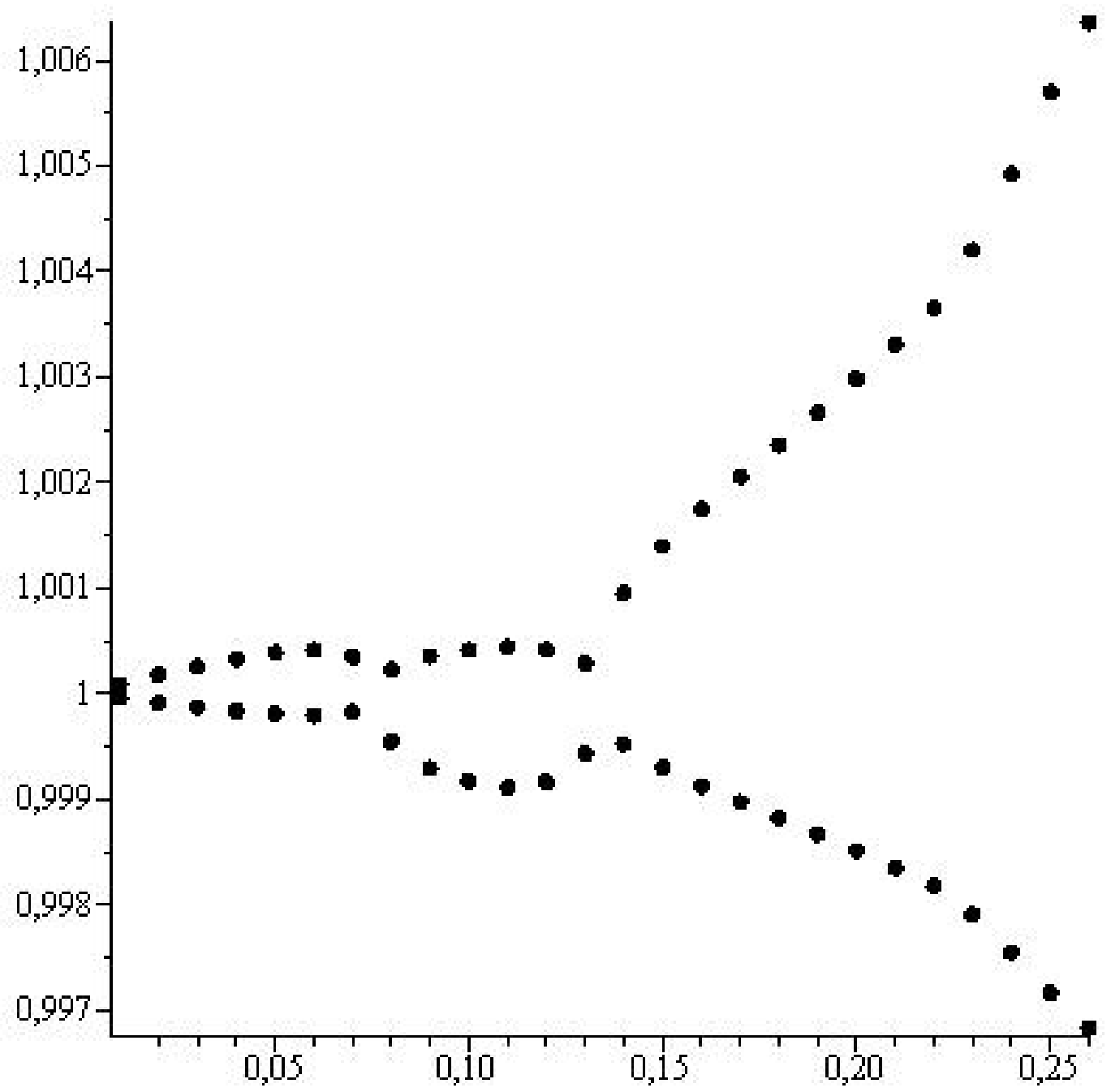}
\end{minipage}
\begin{minipage}{0.45\columnwidth} 
\includegraphics[scale=0.3]{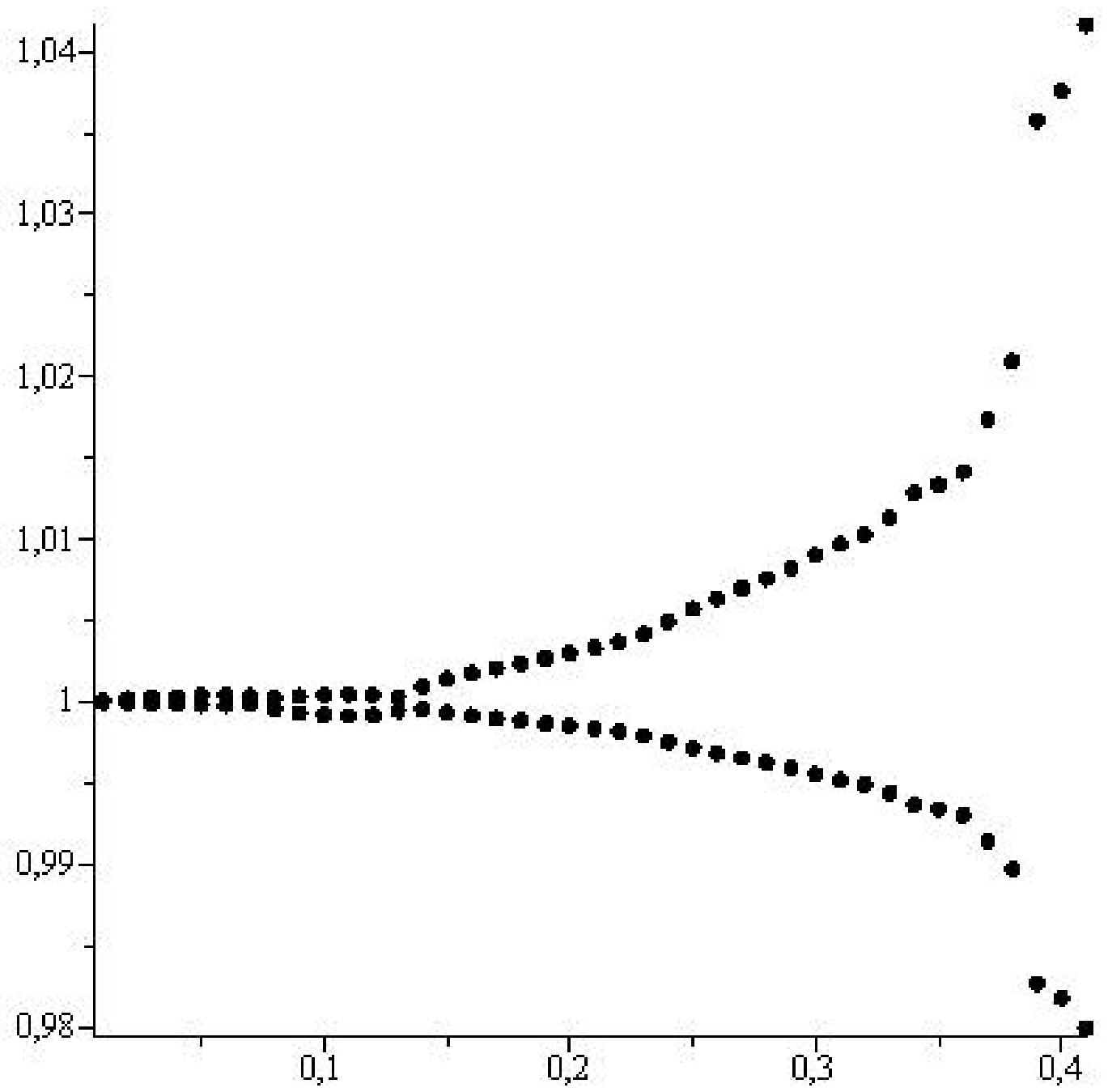}
\end{minipage}
\caption{Values of the characteristic multipliers for the case of non-electrostatic radially symmetric oscillations as a function of $A(0)$, the figure on the left shows the details of the change in the characteristic multipliers with high resolution.}
\label{Pic2}
\end{figure}

Figure \ref{Pic2} shows the dependence  $|\lambda_i(A_*)|$ for different ranges of $A_*$. It is easy to see that the maximum of the eigenvalues exceeds 1 in absolute value, which indicates instability. We see that the quantity $\max\limits_i |\lambda_i(A_*)|-1$, which can be called the measure of instability, changes nonmonotonically on the interval $(0, 0.5)$. However, up to the value of $A_*\approx 0.15$, it is approximately at the same level, being, nevertheless, significantly (two orders of magnitude) larger than the analogous value in the electrostatic case, and then it increases, but not as sharply as in electrostatic case. Among the eigenvalues $\lambda_i$, $i=1, 2, 3$, there are necessarily a pair of complex conjugate ones, and first the real eigenvalue is greater than one in absolute value, then complex conjugate ones is greater in absolute value (for $A_*\in (\approx 0.07, \approx 0.14)$), and then the real eigenvalue again becomes larger in absolute value.

The result on the blow-up of a radially symmetric non-electrostatic perturbation of a general form follows from the same reasoning as in Theorem \ref{T3}.
$\square$
\begin{rem}
The proof that a general perturbation the electrostatic axisymmetric solution in the class of arbitrary affine non-electrostatic solutions (not radially symmetric) collapses in a finite time is carried out in exactly the same way as Theorems \ref{T2} and \ref{T4}, but is more cumbersome, since it requires solving a system of equations of the 9th order and examining 9 eigenvalues. We do not present the results of these calculations.

\begin{rem}
The deviation from the equilibrium  with $\mathcal B_0\ne 0$ behaves quite differently. Indeed,
if in \eqref{10} we replace the representation for $\mathcal B$ with
$\mathcal B(t)=\mathcal B_0+\mathcal B_1(t)\varepsilon^2+o(\varepsilon^2)$,
then we get $A_0(t)=a_0(t)=0$, and the next  terms of expansion are subject to a linear homogeneous system of equations with constant coefficients with a matrix having
purely imaginary  eigenvalues
$
\pm\frac{1}{2} i
\sqrt{4+2\mathcal B_0^2\pm2\sqrt{\mathcal B_0^4+4 \mathcal B_0^2}}$
and zero. Thus, in the first approximation in $\varepsilon$, the solution is a superposition of two periodic motions with different periods. In order to construct the next approximation, one has to solve a linear inhomogeneous equation with constant coefficients. When solving, secular terms arise, but this does not mean that the equilibrium position is unstable (see an example in \cite{Bogolyubov}). Moreover, the numerical results indicate that a small deviations from the equilibrium at $\mathcal B_0\ne 0$ are bounded. However, we do not know an analytical proof of this fact.
 In this case, it is  not possible to apply the method used in the proof of the previous theorems.
 \end{rem}
\begin{figure}[htb!]
\hspace{1.0cm}
\begin{minipage}{0.45\columnwidth}
\includegraphics[scale=0.3]{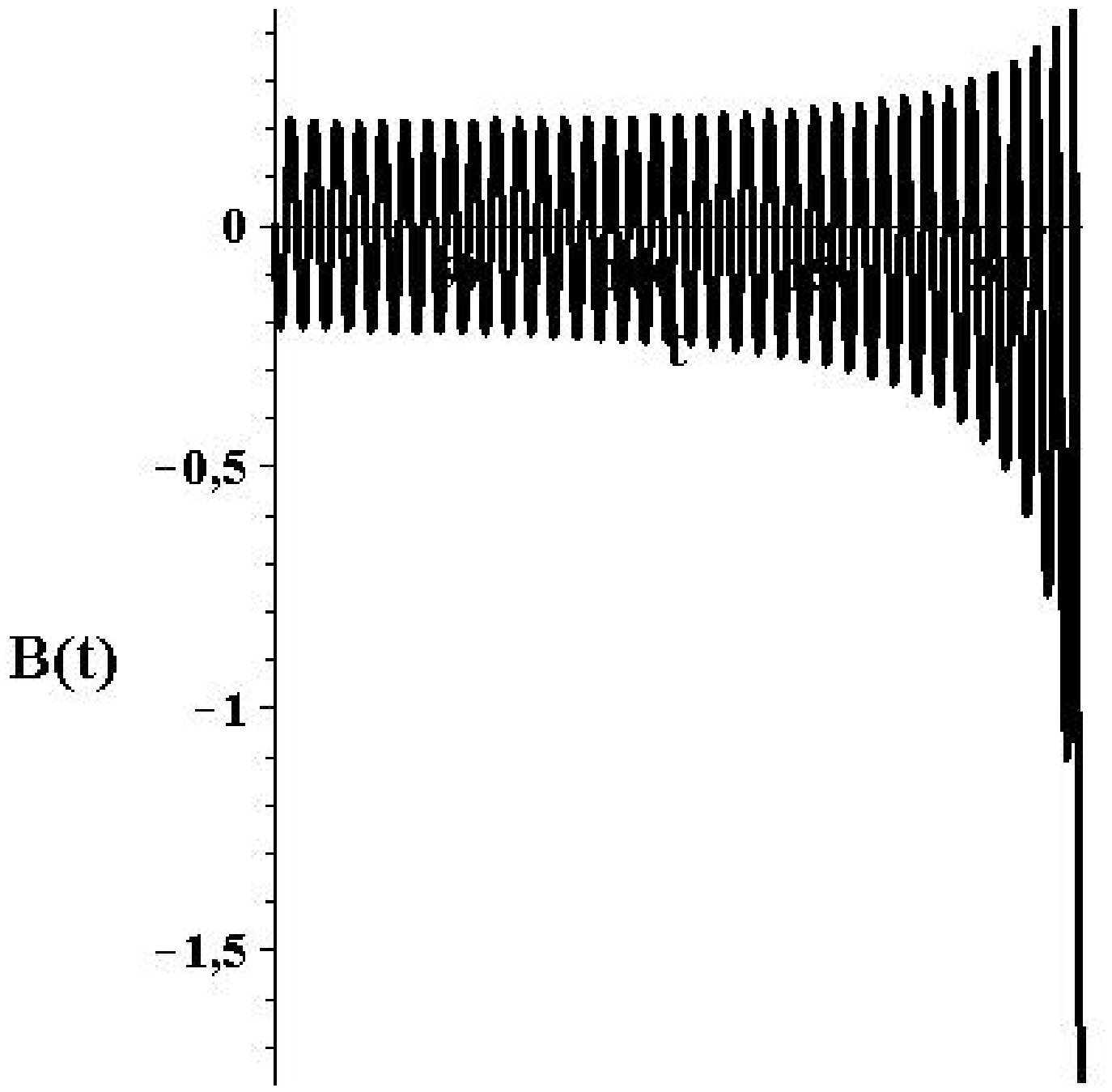}
\end{minipage}
\begin{minipage}{0.45\columnwidth} 
\includegraphics[scale=0.3]{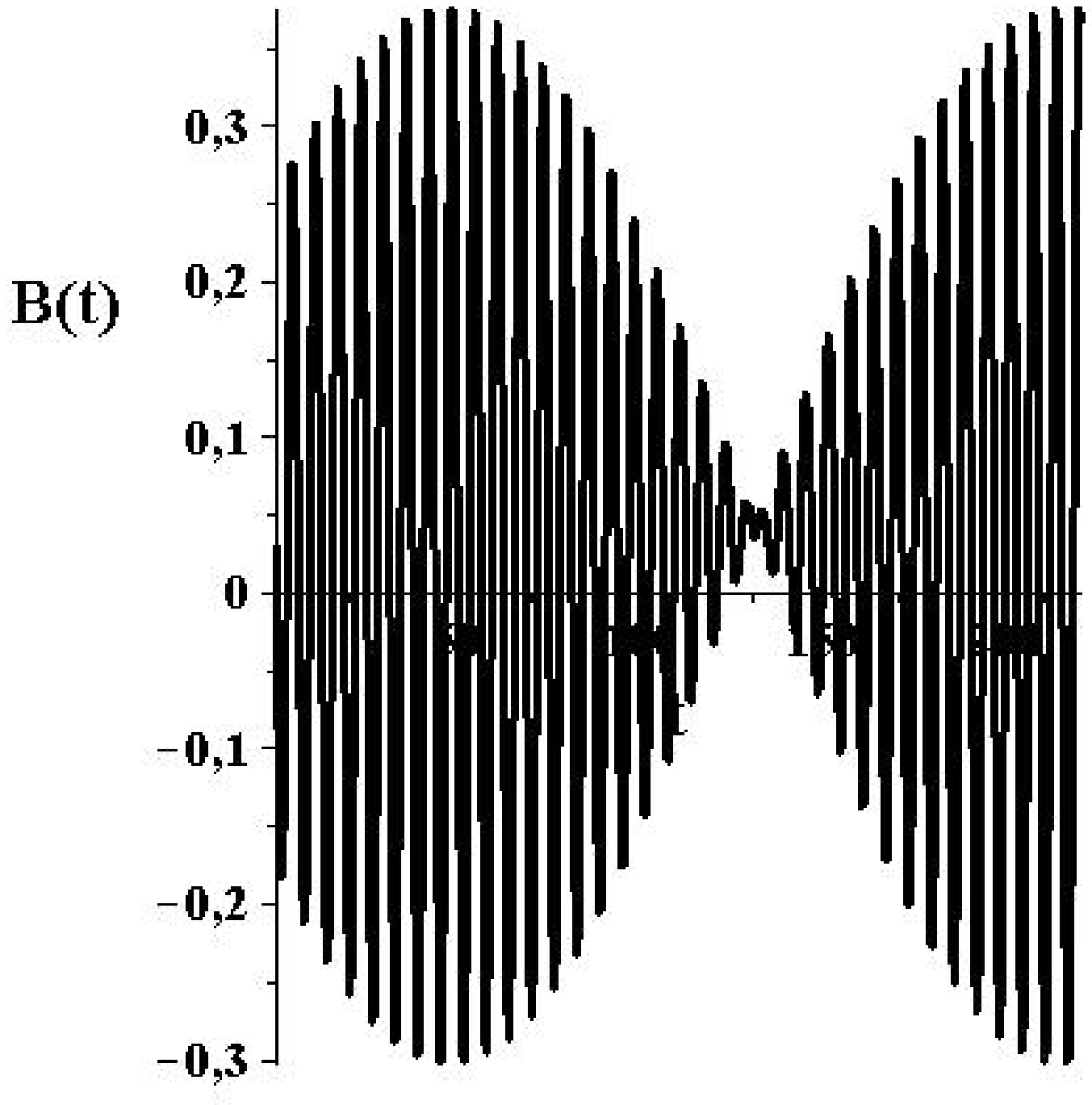}
\end{minipage}
\caption{The difference in the behavior of the solution upon deviation from the equilibrium  at $\mathcal B_0=0$ and $\mathcal B_0\ne 0$.
Initial deviation is chosen as $a(0)=c(0)=0, A(0) =0.1, C(0) = 0.1$, $\mathcal B_0=0$ (left) and $\mathcal B_0= 0.04$
(right). The calculations are done for $t=220$. For $\mathcal B_0= 0$ the solution goes to infinity in finite time.}
\label{Pic3}
\end{figure}
The hypothesis is that the larger $\mathcal B_0$, the wider the neighborhood of the equilibrium, starting from which, the solution remains bounded and globally smooth.

Note that the magnetic field also plays a stabilizing role in other problems related to the description of cold plasma \cite{RCh22}.

 Figure \ref{Pic3} illustrates the difference in the behavior of the magnetic field component for $\mathcal B_0=0$ and $\mathcal B_0 \ne 0$ for the same initial data for the remaining components of the solution.
\end{rem}

\begin{rem} The fact that, for some choice of initial data, the time-dependent coefficients at the second and lower powers of $\varepsilon$ remain bounded for all $t>0$ does not mean that all other coefficients in the expansion of the solution in $\varepsilon$ have the same property. The Floquet theory can be successfully used to prove instability, but it is difficult to apply it to prove stability.
\end{rem}

\section*{Acknowledgments}
Supported by the Moscow Center for
Fundamental and Applied Mathematics under the agreement
№075-15-2019-1621. The authors thanks V.V. Bykov and A.V. Borovskikh
for discussions.

\end{document}